\documentclass[preprint,review,12pt]{elsarticle}
 
\usepackage{graphicx}
\graphicspath{ {./figures/} }
\usepackage{amssymb}
\usepackage{amsthm}
\usepackage{amsmath}
\usepackage[super]{nth}
\usepackage{textcomp}
\usepackage{gensymb}
\usepackage{color, soul}
\usepackage{ulem}
\makeatletter 
\newcount\SOUL@minus
\makeatother  
\usepackage[dvipsnames]{xcolor}
\usepackage{booktabs}
\usepackage{hyperref}
\usepackage{rotating}
\usepackage{pifont}
\usepackage{newunicodechar}
\usepackage{subcaption}
\usepackage{empheq}
\usepackage{dashbox}
\usepackage{lipsum}

\usepackage{listings}
\newcommand{\overbar}[1]{\mkern 1.5mu\overline{\mkern-1.5mu#1\mkern-1.5mu}\mkern 1.5mu}
\usepackage{adjustbox}
\usepackage{makecell}
\usepackage{multirow}
\usepackage{lineno}

\usepackage{float}
\usepackage[most]{tcolorbox}
\usepackage{tikz}
\usepackage{placeins}
\usepackage{booktabs}
\usepackage{tabularx}
\usepackage{amsmath}
\usepackage{hyperref}
\usepackage{setspace}
\usepackage{nomencl}

\makenomenclature
\usepackage{etoolbox}

\newif\iffirstnomgroup
\firstnomgrouptrue

\renewcommand\nomgroup[1]{%
  \iffirstnomgroup
    \global\firstnomgroupfalse
  \else
    \clearpage
  \fi
  \item[\bfseries
  \ifstrequal{#1}{F}{Formula Symbols}{%
  \ifstrequal{#1}{S}{Subscripts}{%
  \ifstrequal{#1}{P}{Superscripts}{}}}%
]}

\newcommand{\nomunit}[1]{%
  \renewcommand{\nomentryend}{\hspace*{\fill}#1}}

\usepackage{tocloft}

\newlength{\mylen}
\renewcommand{\cftfigpresnum}{\figurename\enspace}
\renewcommand{\cftfigaftersnum}{:}
\renewcommand{\cfttabpresnum}{\tablename\enspace}
\newcounter{bullet}

\newcommand{\bulletlabel}[1]{\refstepcounter{bullet}\thebullet\label{#1}}

\renewcommand{\thebullet}{}

\journal{International Journal of Heat and Mass Transfer}

\newlength\myboxwidth

\sethlcolor{SkyBlue}

\begin{document}

\begin{frontmatter}



\title{An Efficient Open-Source Framework for Thermohygromechanical Phase-Field Simulation of Drying Damage Refractory Castables}

\author[DEMA]{M. H. Moreira}\corref{cor1}\ead{murilo.moreira@ufscar.br}
\author[ICMC]{R. F. Ausas}
\author[Grenoble]{S. Dal Pont}
\address[DEMa]{Federal University of Sao Carlos, Materials Microstructural Engineering Group (GEMM), Department of Materials Engineering, Rod. Washington Luiz, km 235, 13565-905, São Carlos, SP Brazil}
\address[Grenoble]{Université Grenoble Alpes, CNRS, Grenoble INP, 3SR, 38000 Grenoble, France}
\address[ICMC]{Institute of Mathematics and Computer Sciences, University of São Paulo, 13566-590, São Carlos, SP, Brazil}
\cortext[cor1]{Corresponding Author}

\begin{abstract}
    The water removal from refractory castables is a complex process involving heat and mass transport and thermohygromechanical stresses. In order to reduce the risk of explosive spalling by improving the heat-up curves, numerical models are critical to predict the castable behavior during drying. However, most of the existing models either neglect the thermomechanical stresses or consider them as an \textit{a posteriori} assessment of the likelihood of failure. The present work proposes an open-source thermohygromechanical phase-field model to simulate the drying of refractory castables and predict the occurrence of explosive spalling. The model is based on an adaptive two-scale time-stepping strategy that allows for a significant reduction of the computational cost associated with the simulation. The model was validated against literature results and qualitative comparisons with experimental data were carried out to assess its ability to predict the failure behavior of partially saturated porous media during heating. The numerical simulations were also applied to investigate the impact of different heating methods on the risk of failure during drying. 
    It was found that even at the significantly lower heating rates typical of refractory castable drying - compared to the fire accidents that motivate most thermomechanical studies on concrete structures - thermomechanical stresses still play a significant role in the explosive spalling, a finding well established for concrete under fire but novel for refractory castables. As a consequence, it was also concluded that the use of heterogeneous heating methods can increase the likelihood of explosive spalling. Ultimately, feasible solutions to increase the drying efficiency and safety that are compatible with the industrial setting were proposed and a free, open-source implementation of the model is provided, enabling both the reproducibility of the present results and its application to optimize drying schedules under different boundary conditions, geometries, and castable compositions.
\end{abstract}

\end{frontmatter}
\clearpage
\printnomenclature
\clearpage
\section{Introduction}
\label{sec:intro}
Refractories are unique structural materials capable of withstanding elevated temperatures, interaction with molten metals, slag, and large mechanical loads \cite{Schacht2004}. Within the different product categories, castables represent a promising alternative to shaped products, thanks to their easier installation and smaller carbon footprint \cite{konigshofer2012developing}).

In the hydraulic bonded castables group, their initial strength after setting is provided by the hydration of its binders \cite{Luz2015book}. However, before firing at high temperatures to develop its final properties, the material must be dried to remove residual water that did not react with the binder, alongside steam generated from the dehydration process of the curing products \cite{Luz2015book}.

This process can last multiple days \cite{Palmer2014Accelerated1} due the use of low heating rates and prolonged periods at fixed temperatures adopted to mitigate the risk of explosive spalling \cite{Palmer2014Accelerated1, luz2021drying, Fey2016b}. This phenomenon is characterized by the sudden removal of large pieces of the refractory lining \cite{Palmer2014Accelerated1, Luz2015book, luz2021drying}. Thus, the heating curves used are often based on semi-empirical knowledge rather than on a rigorous scientific framework \cite{Gong2007}, leading to extended production downtime and higher energy consumption \cite{cunha2023simple}.

In this context, numerical modeling is a viable approach to optimize the heat-up protocols considering a given refractory castable composition, its geometry and dimensions, and the relevant boundary conditions. This potential has motivated several efforts to develop predictive tools for refractory drying over the last decades. Gong \textit{et al.} \cite{Gong1991} proposed a numerical tool to simulate the drying of refractories, based on the framework of Bažant's single-phase thermohygro model, which was initially developed for analyzing unexpected heating of Portland cement concrete (either in fire accidents or nuclear reactors failures) \cite{bazant1978}. A similar strategy was also adopted later by Palmer \textit{et al.} \cite{Palmer2014Accelerated2}.

Fey \textit{et al.} introduced a validated multiphase thermohygro model for refractory castable drying based on earlier developments within civil engineering's studies modeling the high-temperature behavior of concrete structures \cite{Fey2016b}, a common trend observed in the research of refractory drying \cite{moreira2021main}. It was later shown that both Gong's single-phase approach and Fey's multiphase thermohygro models yield similar predictions of the pore pressure developed within refractory castables during drying \cite{moreira2021direct}.

Thus, Cunha \textit{et al.} \cite{cunha2022drying} have used a single phase thermohygro model couples with an \textit{a posteriori} assessment of the likelihood of explosive spalling. The criterion for failure was derived from the comparison of the pressure values with the mechanical strength of the material at specific temperatures within a triaxial stress state calculated from a Weibull analysis. With this model, the effect of different heating rates \cite{cunha2022drying} and a straightforward methodology for heat-up curve selections were proposed \cite{cunha2023simple}. 

These studies enabled valuable insights and guidelines for the drying of refractory castables, however, the mechanical behavior during the drying was either neglected (as in Gong's and Fey's works), or only assumed as an ``a posteriori'' one-way coupling between the thermohygro state and the material's strength. Consequently, the combined effect of moisture clogging (which is the localized reduction of permeability thanks to the condensation of liquid water at colder positions) and the thermomechanical stresses wasn't considered, highlighting the need of a thermohygromechanical model (THM).

Several THM models were introduced in the literature to study the behavior of Portland cement concrete structures at high temperatures \cite{baggio1995thermo, Gawin1999, dal2005experimental, Davie2012a, dal2007multiphase, phan2011finite}. These approaches rely on Biot’s poroelasticity theory to couple the thermohygral and mechanical fields \cite{biot1941general, lewis1999finite}, while differing mainly in their treatment of the mechanical response, ranging from linear elasticity to more sophisticated constitutive descriptions including plasticity, creep and damage \cite{dal2005experimental, dal2007multiphase, Davie2012a, phan2011finite, Gawin2011a}.

Furthermore, the continuum damage model originally proposed by Mazars \textit{et al.} \cite{mazars1986description} is frequently used for predicting the explosive spalling \cite{Gawin1999, dal2007multiphase, Gawin2011, Davie2012a}. This approach assumes an isotropic material behavior, with damage initiating once a strain threshold is surpassed, considering distinct contributions from tension and compression based on the principal strain values \cite{Gawin2011a}. 

The damage is described by a diffuse variable regularized using a non-local formulation based on the averaged equivalent strain over a given volume. This poses challenges in the numerical implementation as it requires access to strain information from neighboring nodes. As highlighted by Arruda \textit{et al.} \cite{arruda2022modified}, this limitation hinders the parallelization of finite-element codes, which is essential for structure-scale simulations. For the same reason, the computations at the boundaries of the domain are not as straightforward as other approaches.

Finally, Cheng \textit{et al.} have proposed the combination of the THM equations with a phase-field model (PFM) for the crack generation and propagation to investigate the fire-induced spalling in Portland cement concrete \cite{cheng2022coupled}. PFM is built upon the Griffiths's theory of brittle fracture of a material with a sharp crack and a regularization of this discontinuity by considering a smeared representation of it depicted by the \textit{phase-field variable} within a region defined by the length scale parameter, $\ell_0$  \cite{francfort1998revisiting, bourdin2000numerical, miehe2010phase, ambati2015review}. Some of the advantages of this method are the possibility of simulating complex fracture processes, such as the crack initiation, propagation, deflection and branching, the capacity of handling complex 3D geometries without requiring any additional \textit{ad-hoc} criteria and a relatively simpler computer implementation in finite element codes \cite{ambati2015review, ghasabeh2023phase}. Additionally, the ease of incorporating it with additional physical process have yield great popularity for PFM in spite of its main drawback, the high computational cost associated with the need for fine meshes the regions surrounding the cracks \cite{ghasabeh2023phase}.

\nomenclature[F]{$\ell_0$}{Length scale parameter \nomunit{[m]}}

Using the proposed THM-PFM, Cheng \textit{et al.} have studied the effect of the length scale parameter, $\ell_0$, the use of different staggered solution algorithms and the effect of the crack on the mass and heat transfer \cite{cheng2022coupled}. It was found out that $\ell_0$ should be twice the minimum element size and the critical energy release rate, $G_c$, should be calibrated to represent both the crack path and the mechanical characteristics properly. In addition, the consideration of the damage effect on the permeability and thermal conductivity led to physically more realistic results. 

\nomenclature[F]{$G_c$}{Critical energy release rate \nomunit{[J/m$^2$]}}

One aspect of Cheng's work that is worth highlighting is the use of fixed time steps for the solution of the THM-PFM equations. This approach led to large computational costs, as the time step size was dictated by the mechanical problem, which is solved in a pseudo-time scale, and not by the thermohygromechanical problem, which is solved in the real time scale. Therefore, even simulations of relatively simple 2D geometries took over three days to predict a period of 300 s \cite{cheng2022coupled}.

Lastly, for the fast heating conditions considered within fire accidents, it was concluded, that the explosive spalling was a result of combined thermal stresses and vapor pressure, with the former playing a major role and the latter accelerating the process.

This last conclusion is in agreement with the findings of previous experimental analysis carried out by Zeiml \cite{zeiml2008experimental} where a high-speed camera was used to estimate the released energy by explosive spalling and THM model studies, as reported by Shen \textit{et al.}, which reported that the primary energy source for the explosive spalling was the thermal stresses, whereas the pore pressure effect was solely a triggering factor \cite{shen2020multiphysics}. Felicetti \textit{et al.} also provided experimental support for this interpretation through fast neutron radiography which allowed the observation the moisture distribution following the depressurization of an internal interface in hot moist concrete, reproducing a crack opening event \cite{felicetti2024fast}. 

Their results show that, upon crack opening, flash vaporization occurs within a thin layer of approximately 1 mm adjacent to the crack surface in as little as 0.03 s, producing temperature drops of 30-90°C at cooling rates up to 120°C/s. Critically, the energy involved in this process, driven by the vast thermal energy stored in the surrouding hot concrete, was found to be orders of magnitude larger than the elastic strain energy available and sufficient to accelerate the fractured fragments to the velocities observed in explosive spalling events. Thus, the pore pressure effect is limited as a trigger of the process.

However, refractory castable drying differs substantially from structures under fire, as the thermal input is controlled and uses significantly lower heating rates. In this slower dynamics it is unclear whether the thermomechanical stresses are still a significant aspect on explosive spalling. 

Furthermore, both the flames in fire accidents and the industrial burners used in castable drying generate intrinsically non-uniform thermal fields due to the complex transport phenomena governing combustion and radiative heat transfer \cite{tien1982flame}. This can be critical both for i) devising boundary conditions that better reproduce the real behavior of concrete structures under fire, because in such conditions the often adopted homogeneous temperature assumption can lead to erroneous conclusions \cite{stern2010experimental} and ii) proposing innovative strategies feasible in industrial scale to increase both safety and efficiency of the refractory castables' drying process, if a reduction of thermal gradient in the direction perpendicular to the drying front is indeed proven to minimize the spalling likelihood.

Therefore, the current work proposes a thermohygromechanical phase-field model for analyzing the drying of refractory castables. The proposed methodology comprises an adaptive two scale time-stepping strategy that reduces the computational cost of the simulation. A benchmark with literature results is presented to validate the model's implementation and qualitative comparisons with experimental results of concrete at high temperatures will be considered to assess the capability of predicting the failure behavior at such conditions. Finally, the model was used to assess the impact of the thermomechanical stresses on the explosive spalling during refractory castables drying in inhomogeneous heating conditions.

\section{Material and Methods}
\label{sec:materialmethods}
A single-phase thermohygromechanical anisotropic phase-field model is proposed for the simulation of the drying of refractory castables. First, the main hypothesis of the model are listed, with the partial derivative equations of the thermohygromechanical and phase-field model presented next. The coupling strategy is detailed and the numerical formulation is provided. The constitutive equations and all the parameters considered in the simulations are listed in the \hyperref[sec:app_const_eq_5CAC]{B}.

\subsection{Physical and Mathematical Model}\label{ssec:mat_model}
\subsubsection{Main Assumptions}\label{sssec:main_assumptions}

The following assumptions are considered in the model proposed herein:
\begin{itemize}
    \item \textbf{A.1}\bulletlabel{A_localeq} The considered phases are in local equilibrium meeting the classical irreversible thermodynamics requirements;
    \item \textbf{A.2}\bulletlabel{A_baz_ads} The mass transport within the partially saturated porous media is dictated by the migration of molecules along the adsorbed layers, which is roughly equivalent of considering the gradient of the water vapor pressure as the driving force of adsorbed water diffusion \cite{bazant1978};
    \item \textbf{A.3}\bulletlabel{A_baz_sp} The transport of liquid and vapor water is represented by the moisture transfer;
    \item \textbf{A.4}\bulletlabel{A_sore_duf} The Soret and the Dufour effect (i.e., the mass flux due to the temperature gradient and the heat flux related to the pressure gradient, respectively) are negligible;
    \item \textbf{A.5}\bulletlabel{A_darcy} Darcy law is assumed to govern the mass transport;
    \item \textbf{A.6}\bulletlabel{A_sorp} The free water (any source of H$_2$O molecules other than the dehydration reactions) content is related to the temperature and pore pressure by the sorption isotherm, and the relationship proposed by Bažant \textit{et al.} is considered \cite{bazant1978};
    \item \textbf{A.7}\bulletlabel{A_perm} The permeability is a function of the temperature and relative humidity (hence, the vapor pressure) as proposed by Bažant \cite{bazant1978}; 
    \item \textbf{A.8}\bulletlabel{A_kinetic_chem} Kinetic effects on the chemical transformations are not considered;
    \item \textbf{A.9}\bulletlabel{A_enthalpy} The enthalpy of desorption is considered equal to the enthalpy of vaporization of liquid water and the dehydration enthalpy is neglected;
    \item \textbf{A.10}\bulletlabel{A_work} The heat related to the work of stresses in concrete is negligible;
    \item \textbf{A.11}\bulletlabel{A_inf_strain} The infinitesimal strain theory is considered for the rigid body;
    \item \textbf{A.12}\bulletlabel{A_elastic} The solid skeleton (the matrix and aggregates that surrounds the pores) is considered as a homogeneous non-deformable, perfectly elastic isotropic material;
    \item \textbf{A.13}\bulletlabel{A_kinetic} Inertial effects are disregarded and the failure mode considered is quasi-static brittle fracture.
\end{itemize}

\subsubsection{Thermohygromechanical Model}\label{sssec:THM_model}
The thermohygro component of the model is based on the single-phase approach developed by Bažant \textit{et al.}, as this simplified approach has shown similar pressure predictions to multiphase models \cite{moreira2021direct}. The mass and heat transfer are given by Equations \ref{eq:mass_transfer_eq} \ref{eq:heat_transfer_eq}, respectively. The full list of the nomenclature considered in the equations is presented at the beginning of this work.

\begin{equation}\label{eq:mass_transfer_eq}
    \frac{\partial{\overbar{\rho}_e}}{\partial{t}} - \nabla \cdot \left( \frac{\kappa}{g} \nabla p_w \right) = \frac{\partial{\overbar{\rho}_d}}{\partial{t}}
\end{equation}

\nomenclature[F]{$\bar{\rho}_i$}{Mass of phase i per unit volume of concrete \nomunit{[kg/m$^3$]}}
\nomenclature[S]{$e$}{Evaporable water}
\nomenclature[S]{$d$}{Dehydration water}
\nomenclature[F]{$t$}{Time \nomunit{[s]}}
\nomenclature[F]{$\kappa$}{Hydraulic conductivity \nomunit{[m/s]}}
\nomenclature[F]{$g$}{Acceleration of gravity \nomunit{[m/s$^2$]}}
\nomenclature[F]{$p_w$}{Pressure of moisture \nomunit{[N/m$^2$]}}

\begin{equation}\label{eq:heat_transfer_eq}
    \rho_{eff} C_{p, eff} \frac{\partial T}{\partial t} - \nabla \cdot (\lambda_{eff} \nabla T) = - C_{p, l}  \frac{\kappa}{g} \nabla p_w \cdot \nabla T + \Delta H_{vap} \frac{\partial \overbar{\rho}_e }{\partial t} 
\end{equation}
\nomenclature[F]{$\rho$}{Density \nomunit{[kg/m$^3$]}}
\nomenclature[S]{$eff$}{Effective}
\nomenclature[F]{$C_{p, i}$}{Isobaric specific heat of phase i \nomunit{[J/(kg K)]}}
\nomenclature[F]{$T$}{Temperature \nomunit{[K]}}
\nomenclature[F]{$\lambda_{eff}$}{Effective thermal conductivity \nomunit{[W/(m K)]}}
\nomenclature[S]{$l$}{Liquid water}
\nomenclature[F]{$\Delta H_{vap}$}{Enthalpy of evaporation of water \nomunit{[J/kg]}}

\nomenclature[F]{$\epsilon$}{Emissivity \nomunit{[-]}}
\nomenclature[F]{$\sigma$}{Stefan-Boltzmann's constant \nomunit{[5.67 10$^{-8}$ W/(m$^2$ K$^4$)]}}
\nomenclature[F]{$w_i$}{Initial water content \nomunit{[kg/m$^3$]}}
\nomenclature[S]{$env$}{Environment}
\nomenclature[F]{$S_l$}{Pore saturation with liquid \nomunit{[-]}}
\nomenclature[F]{$h_T$}{Heat transfer coefficient \nomunit{[W/m$^2$K]}}
\nomenclature[F]{$h_m$}{Mass transfer coefficient \nomunit{[s/m]}}
\nomenclature[F]{$\Delta m_{dehyd}$}{Mass released by dehydration \nomunit{[kg/m$^3$]}}

\subsubsection{Phase-Field Model}\label{sssec:PFM}
Starting with the poroelasticity theory proposed by Biot  \cite{biot1941general} - with the effective stress represented by Equation \ref{eq:effective_stress} - and combining it with the thermomechanical stresses developed inside the material - shown in Equation \ref{eq:thermomechanical_stress} - the effective strain,  $\boldsymbol{\varepsilon}_{\textrm{eff}}$, resulting from the contributions of both the fluids exerting pressure $p_f$ over the pore internal walls and the differences of thermal expansion are calculated by Equation \ref{eq:effective_strain_decomp}.
\begin{equation}\label{eq:effective_stress}
    \boldsymbol{\sigma}_{\textrm{eff}} = \boldsymbol{\sigma}_{\textrm{el}} + \mathbf{I} p_f
\end{equation}
where $\boldsymbol{\sigma}_{\textrm{eff}}$ is the effective stress resulting from elastic stresses, $\boldsymbol{\sigma}_{\textrm{el}}$ and the fluids exerting a given pressure $p_f$ within the pore walls inside the concrete. $\mathbf{I}$ is the identity tensor.
\nomenclature[F]{$\boldsymbol{\sigma}$}{Stress \nomunit{[Pa]}}
\nomenclature[S]{$el$}{Elastic}
\nomenclature[F]{$\mathbf{I}$}{Identity tensor \nomunit{[-]}}
\nomenclature[F]{$p_f$}{Fluid pressure \nomunit{[Pa]}}

\begin{equation}
    \label{eq:thermomechanical_stress}
    \boldsymbol{\sigma}_{\textrm{th}} = \frac{K}{3} \alpha \Delta T \boldsymbol{I}
\end{equation}
with $K$ being the bulk modulus, $\alpha$ the thermal expansion, $\Delta T$ the difference in temperature in comparison to a reference state and $\boldsymbol{I}$ the identity tensor.
\nomenclature[S]{$th$}{Thermomechanical}
\nomenclature[F]{$K$}{Bulk modulus \nomunit{[Pa]}}
\nomenclature[F]{$\alpha$}{Thermal expansion \nomunit{[1/K]}}

The effective strain is decomposed into three contributions, as shown in Equation \ref{eq:effective_strain_decomp}.
\begin{equation}
    \label{eq:effective_strain_decomp}
    \boldsymbol{\varepsilon}_{\textrm{eff}} = \boldsymbol{\varepsilon}_{\textrm{el}} + \boldsymbol{\varepsilon}_{\textrm{hyd}} + \boldsymbol{\varepsilon}_{\textrm{th}}
\end{equation}
where $\boldsymbol{\varepsilon}_{\textrm{el}}$ is the elastic deformation strain due to mechanical loads over the representative volume element  of the refractory castable, $\boldsymbol{\varepsilon}_{\textrm{hyd}}$ is the ``hydraulic strain'', which represents the strain by the vapor pressurization within the pores and is given by Equation \ref{eq:hydraulic_strain}
\begin{equation}\label{eq:hydraulic_strain}
    \boldsymbol{\varepsilon}_{\textrm{hyd}} = \frac{b}{K} \ p_f \ \mathbf{I} 
\end{equation}
where $b$ is the Biot coefficient, the bulk modulus is $K$ and $p_f$ the fluid pressure. 
\nomenclature[F]{$\boldsymbol{\varepsilon}$}{Strain \nomunit{[-]}}
\nomenclature[S]{$hyd$}{Hydraulic}
\nomenclature[F]{$b$}{Biot coefficient \nomunit{[-]}}
The Biot poroelasticity coefficient, which describes the extent that the fluid pressure causes strains on the solid skeleton, is the ratio of the fluid and bulk solid volume changes, considering a fixed pore pressure \cite{raziperchikolaee2020effect}. This last assumption only holds in the case currently studied because of the coexistence of two different timescales. 
The thermohygro behavior is transient, and the involved phenomena develop in the real timescale, whereas hypothesis \hyperref[A_kinetic]{A.13} means that the mechanical problem is solved as a stationary one, in a pseudo-time scale [see Figure \ref{fig:coupling} b)], with the mechanical loading resulting from the thermomechanical and hydraulic process imposed gradually over the domain. Therefore, in the mechanical model, the pore pressure is considered constant between successive physical time steps.
Finally, $\boldsymbol{\varepsilon}_{\textrm{th}}$ is the strain resulting from the thermal expansion of the solid given by Equation \ref{eq:thermomechanical_strain}.
\begin{equation}\label{eq:thermomechanical_strain}
    \boldsymbol{\varepsilon}_{\textrm{th}} = \alpha \ \Delta T \ \mathbf{I} 
\end{equation}
with the assumptions of an isotropic material (\hyperref[A_elastic]{A.12}) and linear thermal expansion coefficient represented by $\alpha$ whereas $\Delta T$ is the difference in temperature in comparison to a reference state.

The strain resulting from the combined effects of mechanical loads, thermomechanical stresses, and hydraulic strains, $\boldsymbol{\varepsilon}_{\textrm{eff}}$, will cause a variation in the domain's potential elastic energy, denoted as $E_{\textrm{el}}$. Considering all other energy contributions for a specific domain, the overall energy balance for the system can be expressed as:
\begin{equation}\label{eq:gen_energy_balance}
    E_\textrm{T} = E_{\textrm{el}} + E_{\textrm{surf}} + E_{\textrm{kin}} + E_{\textrm{dis}} + W_{\textrm{ext}}
\end{equation}
In this equation, $E_\textrm{T}$ represents the total energy, $E_{\textrm{surf}}$ is the energy associated with the crack surface, $E_{\textrm{kin}}$ accounts for the kinetic energy, $E_{\textrm{dis}}$ represents the dissipated energy in the form of sound, light, or plastic deformation, and $W_{\textrm{ext}}$ represents the work done by external forces.

\nomenclature[F]{$E$}{Energy \nomunit{[J]}}
\nomenclature[S]{$T$}{Total}
\nomenclature[S]{$\textrm{surf}$}{Surface}
\nomenclature[S]{$\textrm{kin}$}{Kinetic}
\nomenclature[S]{$\textrm{dis}$}{Dissipated}
\nomenclature[F]{$W$}{Work \nomunit{[J]}}
\nomenclature[S]{$\textrm{ext}$}{External}

Considering assumptions \hyperref[A_elastic]{A.12} and \hyperref[A_kinetic]{A.13}, the dissipated and kinetic energies are neglected ($E_{\textrm{dis}} = 0$ and $E_{\textrm{kin}}=0$). Furthermore, in the staggered approach adopted herein, the mechanical problem takes into account the fluid pressure and thermal expansion effects as external forces [see Figure \ref{fig:coupling} a)] \cite{cheng2022coupled}.

Thus, in the THM-PF developed herein, the energy functional is given by:
\begin{equation}\label{eq:energy_balance}
    E_{\textrm{tot}}(\boldsymbol{u}, \Gamma_c)= \underbrace{\int_{\Omega} \psi_{\textrm{el}}(\boldsymbol{\varepsilon}_{\textrm{eff}}(\boldsymbol{u})) \mathrm{d}\Omega}_{E_{\textrm{el}}} + \underbrace{\int_{\Gamma_c} G_c \mathrm{d} \Gamma}_{E_{\textrm{surf}}} + \underbrace{\int_{\Omega} \boldsymbol{b}_t \cdot \boldsymbol{u} \mathrm{d} \Omega-\int_{\partial \Omega_N} \boldsymbol{f}_t \cdot \boldsymbol{u} \mathrm{d} S}_{W_{\textrm{ext}}}
\end{equation}
where $\Omega$ represents the domain of interest, $\psi_{\textrm{el}}(\boldsymbol{\varepsilon}(\boldsymbol{u}))$ denotes the elastic strain energy density, $\Gamma_c$ represents the crack surface, $G_c$ represents Griffith's critical energy release rate, $\boldsymbol{b}_t$ indicates the body forces (such as forces due to gravity, electric or magnetic fields), and $\boldsymbol{f}_t$ represents surface forces.

\nomenclature[F]{$\boldsymbol{u}$}{Displacement vector \nomunit{[m]}}
\nomenclature[F]{$\Gamma_c$}{Crack surface \nomunit{[m$^2$]}}
\nomenclature[F]{$\psi$}{Elastic strain energy density \nomunit{[J/m$^3$]}}
\nomenclature[F]{$\Omega$}{Domain of interest \nomunit{[m$^3$]}}
\nomenclature[F]{$S$}{Boundary of the domain \nomunit{[m$^2$]}}
\nomenclature[F]{$\boldsymbol{b}_t$}{Body forces \nomunit{[N/m$^3$]}}
\nomenclature[F]{$\boldsymbol{f}_t$}{Surface forces \nomunit{[N/m$^2$]}}

The linear momentum balance equation with the displacement vector as the primary variable and the Neumann boundary conditions defined over $\Gamma_N^{\textrm{PFM}}$ is presented in Equation \ref{eq:lin_momentum}:
\begin{equation}\label{eq:lin_momentum}
\begin{cases}
\nabla \cdot \boldsymbol{\sigma}_{\textrm{el}}(\boldsymbol{u}) + \mathbf{b}_t = 0 & \text { in } \Omega \\ \boldsymbol{\sigma}_{\textrm{el}}(\boldsymbol{u})  \cdot \boldsymbol{\vec{n}} = \mathbf{f}_t & \text { on } \Gamma_N^{\textrm{PFM}}
\end{cases}   
\end{equation}
with $\boldsymbol{\vec{n}}$ being the unit normal vector.

\nomenclature[F]{$\boldsymbol{\vec{n}}$}{Unit normal vector \nomunit{[-]}}
\nomenclature[P]{$\textrm{PFM}$}{Relative to the phase-field model}
\nomenclature[P]{$\textrm{TH}$}{Relative to the thermohygro model}

\nomenclature[S]{$N$}{Neumann boundary conditions}
\nomenclature[S]{$D$}{Dirichlet boundary conditions}
\nomenclature[S]{$R$}{Robin boundary conditions}

An anisotropic phase-field approach is adopted where the stress tensor is divided into the crack driving component (i.e. the tensile one), from the compressive ones. 

Hence, the elastic stresses are calculated as follows:
\begin{equation}\label{eq:sigma_el}
    \boldsymbol{\sigma}_{\textrm{el}}(\boldsymbol{\varepsilon}, \phi) = \frac{\partial \psi(\boldsymbol{\varepsilon}_{\textrm{el}}, \phi)}{\partial \boldsymbol{\varepsilon}_{\textrm{el}}} = g_d(\phi) \frac{\partial \psi_D(\boldsymbol{\varepsilon}_{\textrm{el}})}{\partial \boldsymbol{\varepsilon}_{\textrm{el}}} + \frac{\partial \psi_R(\boldsymbol{\varepsilon}_{\textrm{el}})}{\partial \boldsymbol{\varepsilon}_{\textrm{el}}}
\end{equation}
with $\phi$ being the phase-field variable representing the damage, varying from 0 (intact) to 1 (fully damaged), $g$ representing the reduction of the stiffness as a function of the damage, and $\psi_D$ and $\psi_R$ are the degradable and residual components of the elastic strain energy
density.

\nomenclature[F]{$\phi$}{Phase-field variable representing damage \nomunit{[-]}}
\nomenclature[F]{$g_d$}{Degradation function \nomunit{[-]}}
\nomenclature[S]{$D$}{Degradable}
\nomenclature[S]{$R$}{Residual}

The decomposition of the elastic strain energy density into its degradable and residual components is accomplished by the spectral decomposition method, resulting into Equations \ref{eq:psi_D} and \ref{eq:psi_R}
\begin{equation}\label{eq:psi_D}
    \psi_D(\boldsymbol{\varepsilon})=\frac{1}{2} \lambda\langle\operatorname{tr}(\boldsymbol{\varepsilon})\rangle_{+}^2+\mu \boldsymbol{\varepsilon}^{+} \cdot \boldsymbol{\varepsilon}^{+}
\end{equation}
\begin{equation}\label{eq:psi_R}
    \psi_R(\boldsymbol{\varepsilon})=\frac{1}{2} \lambda\langle\operatorname{tr}(\boldsymbol{\varepsilon})\rangle_{-}^2+\mu \boldsymbol{\varepsilon}^{-} \cdot \boldsymbol{\varepsilon}^{-}
\end{equation}
where $\lambda$ and $\mu$ are Lamé's first and second constants, $\operatorname{tr}(\textbf{\textrm{T}})$ is the trace of a given tensor $\textbf{\textrm{T}}$ and the operators $\langle \ \ \rangle_{+}$ and $\langle \ \ \rangle_{-}$ represent the positive and negative parts of a given scalar, vector or tensor field, following the Macauley bracket notation.
\nomenclature[F]{$\lambda$}{Lamé's first constant \nomunit{[Pa]}}
\nomenclature[F]{$\mu$}{Lamé's second constant \nomunit{[Pa]}}
In Equation \ref{eq:psi_D}, the tensor $\boldsymbol{\varepsilon}^{+}$ denotes the positive component of the strain tensor. It exclusively considers tensile strains (positive principal strains) and disregards compressive ones (negative principal strains). Therefore, $\boldsymbol{\varepsilon}^{+}$ is obtained by setting all the negative eigenvalues (principal strains) of $\boldsymbol{\varepsilon}$ to zero and retaining only the positive ones. This can be expressed as follows:
\begin{equation}\label{eigens_p}
    \boldsymbol{\varepsilon}^{+}=\sum_i\left\langle e_i\right\rangle_{+} \boldsymbol{e}_i \otimes \boldsymbol{e}_i
\end{equation}
In this equation, $\mathbf{e}_i$ and $e_i$ represent the eigenvectors and corresponding eigenvalues of the strain tensor, respectively.
Conversely, $\boldsymbol{\varepsilon}^{-}$ on Equation \ref{eq:psi_R}, defined as
\begin{equation}\label{eigens_n}
\boldsymbol{\varepsilon}^{-}=\sum_i\left\langle e_i\right\rangle_{-} \boldsymbol{e}_i \otimes \boldsymbol{e}_i
\end{equation}
considers only the negative part of the strain tensor.

\renewcommand\arraystretch{1}

The degradation function, $g_d(\phi)$, is assumed to behave parabolic, following the AT$_\textrm{2}$ approach \cite{Vicentini2024}:
\begin{equation}\label{eq:degrad}
    g_d(\phi) = (1 - \phi)^2 + \varsigma
\end{equation}
with $\varsigma$ being a numeric residual stiffness adopted to be $\varsigma = 10^{-6}$ in the current work to avoid convergence issues when $\phi \to 1$.

\nomenclature[F]{$\varsigma$}{Numeric residual stiffness \nomunit{[-]}}

Based on the regularized approach by Bourdin et al. \cite{bourdin2000numerical}, the crack surface energy, $E_{\text{surf}}$, is given by:
\begin{equation}\label{eq:E_surf}
    E_{\textrm{surf}}(\phi) = \int_{\Gamma_c} G_c \mathrm{d} \Gamma = G_c \Gamma(\phi) = \int_{\Omega}\left(\frac{G_c}{2 \ell_0} \phi^2+\frac{\ell_0}{2}|\nabla \phi|^2\right) d \Omega
\end{equation}
where $\ell_0$ is the characteristic length.

Considering the minimization of the resulting total energy functional, $E_{\textrm{tot}} = E_{\textrm{el}} + E_{\textrm{surf}}$, and introducing the history variable, $\mathcal{H}$ to account for the crack irreversibility (as defined in Equation \ref{eq:history}), the phase-field problem is given by:
\begin{equation}\label{eq:phi_diffusion}
\begin{dcases}
\frac{G_c}{\ell_0} (\phi - \ell_0^2 \nabla^2\phi) = 2 (1 - \phi) \mathcal{H} & \text { in } \Omega \\ 
\nabla \phi \cdot \boldsymbol{\vec{n}} = 0 \vphantom{\frac{0}{0}} & \text { on } \Gamma^{PFM}
\end{dcases}   
\end{equation}
\begin{equation}\label{eq:history}
    \mathcal{H} = \max_{t \in [0, \tau]}\psi_D(\boldsymbol{\varepsilon}(t))
\end{equation}

\nomenclature[F]{$\mathcal{H}$}{History variable for crack irreversibility \nomunit{[J/m$^3$]}}

\definecolor{orangered}{RGB}{255,69,0}
\definecolor{navy}{RGB}{0,0,128}
\definecolor{dodgerblue}{RGB}{30, 144, 255}
\definecolor{goldenrod}{RGB}{218, 165, 32}

Hence, the resulting thermohygromechanical phase-field model comprises 4 primary variables, the pore water vapor \textcolor{dodgerblue}{$p_w$}, the temperature \textcolor{orangered}{$T$}, the displacement vector \textcolor{goldenrod}{$\boldsymbol{u}$} and the phase-field variable indicating the damage \textcolor{navy}{$\phi$}, and 4 equations:
\begin{align}\label{eq:summary_strong}
    \frac{\partial{\overbar{\rho}_e}}{\partial{t}}  - \nabla \cdot \left( \frac{\kappa}{g} \nabla\textcolor{dodgerblue}{p_w} \right) -\frac{\partial{\overbar{\rho}_d}}{\partial{t}} &= 0 \phantom{ab} \text { in } \Omega \\
    \rho_{\textrm{eff}} C_{p, \textrm{eff}} \frac{\partial \textcolor{orangered}{T}}{\partial t} -\nabla \cdot (\lambda_{\textrm{eff}} \nabla \textcolor{orangered}{T}) + C_{p, l}  \frac{\kappa}{g} \nabla \textcolor{dodgerblue}{p_w} \cdot \nabla \textcolor{orangered}{T} - \Delta H_{\textrm{vap}} \frac{\partial \overbar{\rho}_e }{\partial t} &= 0 \phantom{ab} \text { in } \Omega \\
    \nabla \cdot \boldsymbol{\sigma}_{\textrm{el}} (\textcolor{goldenrod}{\boldsymbol{u}}, \textcolor{navy}{\phi}) + \mathbf{b}_t &= 0 \phantom{ab} \text { in } \Omega \\
    \frac{G_c}{\ell_0} (\textcolor{navy}{\phi} - \ell_0^2 \nabla^2\textcolor{navy}{\phi}) - 2 (1 - \textcolor{navy}{\phi}) \mathcal{H}(\textcolor{goldenrod}{\boldsymbol{u}}) &= 0 \phantom{ab} \text { in } \Omega
\end{align}

\subsection{Coupling}\label{sssec:coupling}
\begin{figure}[H]
\centering
\includegraphics[width=1\linewidth]{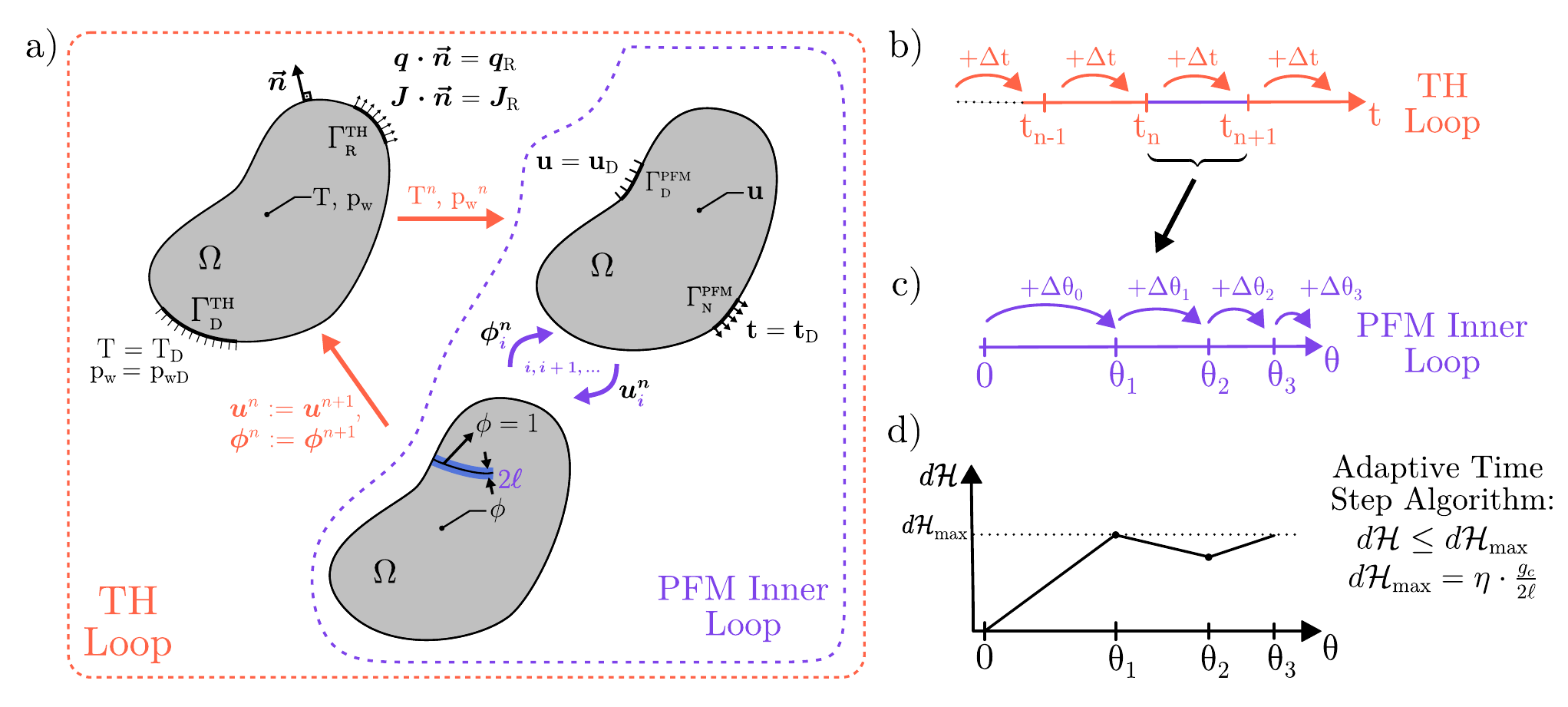}
\caption{a) Schematic representation of the coupled THM-PFM domain and the respective boundary conditions. \label{fig:coupling}}
\end{figure}

The thermohygro equations are solved using a monolithic fully coupled strategy, whereas a staggered approach is considered for the displacement and phase-field problems, as shown in Figure \ref{fig:coupling} a). The domain $\Omega$ is considered, where two scalar fields exist, one related to the temperature distribution, $T$, and the other to the moisture pressure, $p_w$. 
The mass and heat exchange by natural convection with the environment are modeled as Robin type boundary conditions on a part of 
$\partial\Omega$ ($=\Gamma$) denoted by $\Gamma^{\textrm{TH}}_{R}$. Dirichlet boundary conditions are applied over $\Gamma^{\textrm{TH}}_{D}$. The remainder of the boundary $\Gamma^{\textrm{TH}}_{N,0}$ is subject to homogeneous Neumann boundary condition, such that $\Gamma = \Gamma^{\textrm{TH}}_{R} \cup \Gamma^{\textrm{TH}}_{D} \cup \Gamma^{\textrm{TH}}_{N,0}$. For the mechanical problem, displacement boundary conditions are enforced on $\Gamma^{\textrm{PFM}}_{D}$, whereas on $\Gamma^{\textrm{PFM}}_{N}$, a traction $\mathbf{t}$ is naturally imposed. On the rest of $\Gamma$ the natural traction free condition holds. 

Finally, the phase-field equations are considered for the representation of the diffuse crack, where the fully damaged region is taken as $\phi=1$, and $\phi=0$ in the intact parts. A homogeneous Neumann boundary condition is considered for the phase-field problem. Note that the problem is still well-posed due to the zero-th order term
$\frac{G_c}{\ell_0} \phi$ present in Equation (\ref{eq:phi_diffusion}).

Given the assumption of quasi-static brittle fracture and neglecting inertial effects (\hyperref[A_kinetic]{A.13}), two different timescale are considered, one related to the real time of the transient heat and mass transfer and a second one, regarded as a pseudo-time scale which dictates the rate of application of the displacement resulting from the pressure and thermomechanical effects. This new strategy is not adopted by Cheng et al., where a constant step - which was constrained to small values, such as $10^{-2}$ s - was considered. This resulted in large computational costs for even relatively simple geometries and boundary conditions \cite{cheng2022coupled}.

In the simulation of refractory drying, a significant part of the simulated process is related to the initial stages of drying. In this stage, thermal gradients and vapor pressure induce only minor displacements, and small increases on the stored elastic energy density, allowing for higher rates of application of the thermo and hygro displacements.

To have an even more efficient approach, one should consider that in the specific case, a large portion of the simulation concerns the initial stages of heating, where both the thermal gradients and the vapor pressure results in small displacements. This scenario means that if a constant and sufficiently small pseudo-time step $\Delta \theta$ is adopted to guarantee that the numerical results are accurate during the crack initiation and propagation, a large portion of the total simulation run time is wasted as larger time steps could be safely used in the earlier moments of the heating.

Thus, an adaptive pseudo-time step strategy can be used, such as the one proposed by Gupta et al. and Molnar et al. \cite{gupta2020auto, molnar2020toughness}, where the increment of elastic strain energy is controlled as described in Equation \ref{eq:AAPTS}.
\begin{equation}\label{eq:AAPTS}
    \textrm{d}\mathcal{H} \leq \ \eta \frac{G_c}{2 \ell_0}
\end{equation}
with $\eta$ being a scalar assumed to be $\eta = 0.5$.

\nomenclature[F]{$\eta$}{Scalar controlling the increment of elastic strain energy in the adaptive pseudo-time step strategy}

\subsection{Numerical Formulation}\label{ssec:num_form}
The finite element method is considered for the spatial discretization of the domain, following the Galerkin approach. The variational formulation is presented in Equations \ref{eq:summary_weak}: 
Find fields $T$, $p_w$, $\boldsymbol{u}$ and $\phi$,
such that
\begin{align}\label{eq:summary_weak}
    \begin{split}
    \int_{\Omega} \left[\left( \frac{\partial{\overbar{\rho}_e}}{\partial{t}} - \frac{\partial{\overbar{\rho}_d}}{\partial{t}} \right) v_1 + \frac{\kappa}{g} \nabla \textcolor{dodgerblue}{p_w} \cdot \nabla v_1 \right] \, \textrm{d}\Omega \ &- \\
    - \int_{\Gamma} \left( \frac{\kappa}{g} \nabla \textcolor{dodgerblue}{p_w} \cdot \boldsymbol{\vec{n}}\right) \ v_1 \, \textrm{d}\Gamma &= 0
    \end{split}
    \\
    \begin{split}
    \int_{\Omega} \left( \rho_{\textrm{eff}} C_{p, \textrm{eff}} \frac{\partial \textcolor{orangered}{T}}{\partial t}
    + C_{p, l} \frac{\kappa}{g} \nabla \textcolor{dodgerblue}{p_w} \cdot \nabla \textcolor{orangered}{T} 
    - \Delta H_{\textrm{vap}} \frac{\partial \overbar{\rho}_e }{\partial t} \right) v_2 \ &+ \\ + \int_{\Omega} \lambda_{\textrm{eff}} \nabla \textcolor{orangered}{T} \cdot \nabla v_2 \, \textrm{d}\Omega \ - \int_{\partial \Omega} \left( \lambda_{\textrm{eff}} \nabla \textcolor{orangered}{T} \cdot \boldsymbol{\vec{n}}\right) \ v_2 \, \textrm{d}\Gamma &= 0
    \end{split}
     \\
    \int_{\Omega} \left( \boldsymbol{\sigma}_{\textrm{el}} (\textcolor{goldenrod}{\boldsymbol{u}}, \textcolor{navy}{\textcolor{navy}{\phi}}) : \nabla \boldsymbol{v}_3 
    - \mathbf{b}_t \cdot \boldsymbol{v}_3 \right) \textrm{d}\Omega 
    - \int_{\partial \Omega} \boldsymbol{\sigma}_{\textrm{el}} (\textcolor{goldenrod}{\boldsymbol{u}}, \textcolor{navy}{\phi}) \cdot \boldsymbol{\vec{n}} \cdot \boldsymbol{v}_3 \, \textrm{d}\Gamma &= 0 \\
    \int_{\Omega} \left( \frac{G_c}{\ell_0} \textcolor{navy}{\phi} 
    - 2 (1 - \textcolor{navy}{\phi}) \mathcal{H}(\textcolor{goldenrod}{\boldsymbol{u}}) \right) v_4 + G_c \ell_0 \nabla \textcolor{navy}{\phi} \cdot \nabla v_4 \, \textrm{d}\Omega  &= 0
\end{align}
$\forall v_1$, $v_2$, $\boldsymbol{v}_3$ and $v_4$. Above, $\boldsymbol{\vec{n}}$ denote the unit exterior normal to $\Gamma$.
\nomenclature[F]{$v_i$}{Test function for the $i$-th primary variable \nomunit{[-]}}
As previously mentioned, in the permeable parts of the boundary, the mass transfer is described considering a Robin type boundary condition, as defined in Equation \ref{eq:BC_perm}.
\begin{equation}\label{eq:BC_perm}
    - \frac{\kappa}{g} \nabla \textcolor{dodgerblue}{p_w} \cdot \boldsymbol{\vec{n}} = \beta_\textrm{M} (\textcolor{dodgerblue}{p_w} - p_{w, \infty})  \phantom{a} \text { on } \Gamma^{\textrm{TH}}_{\textrm{R}, \textrm{M}}
\end{equation}
where the mass transfer coefficient with units of [s/m] is  $\beta_\textrm{M}$, and $p_{w, \infty}$ represents the water vapor pressure at external environment.
\nomenclature[F]{$\beta_\textrm{M}$}{Mass transfer coefficient \nomunit{[s/m]}}
\nomenclature[F]{$\beta_\textrm{T}$}{Heat transfer coefficient \nomunit{[W/(m$^2$ K)]}}
\nomenclature[S]{$\infty$}{Symbol representing the external environment}
Whereas, the cooling by natural convection is represented by:
\begin{equation}\label{eq:BC_conv}
    - \lambda \nabla \textcolor{orangered}{T} \cdot \boldsymbol{\vec{n}} = \beta_\textrm{T} (\textcolor{orangered}{T} - T_{\infty})  \phantom{a} \text { on } \Gamma^{\textrm{TH}}_{\textrm{R}, \textrm{T}}
\end{equation}
with $\beta_\textrm{T}$ being the heat transfer coefficient, and $T_{\infty}$ the temperature at the surroundings.

Equations \ref{eq:BC_impermeable} and \ref{eq:BC_adiabatic} represent the portions of the boundary where homogeneous Neumman boundary conditions are assumed, 

\begin{equation}\label{eq:BC_impermeable}
    - \frac{\kappa}{g} \nabla \textcolor{dodgerblue}{p_w} \cdot \boldsymbol{\vec{n}} = 0  \phantom{a} \text { on } \Gamma^{\textrm{TH}}_{\textrm{NH}, \textrm{M}}
\end{equation}
\begin{equation}\label{eq:BC_adiabatic}
    - \lambda \nabla \textcolor{orangered}{T} \cdot \boldsymbol{\vec{n}} = 0 \phantom{a} \text { on } \Gamma^{\textrm{TH}}_{\textrm{NH}, \textrm{T}}
\end{equation}

\nomenclature[S]{$NH$}{Homogeneous Neumann boundary condition}

Note that the test functions $v_1$, $v_2$ and the appropriate components of $\boldsymbol{v}_4$ are assumed to vanish on the portions of $\partial\Omega$ where Dirichlet boundary conditions are imposed on $p_w$, $T$ and $\boldsymbol{u}$, respectively. In all numerical experiments in this work, standard $P_1$ (Lagrangian) elements are considered.
Finally, the FEniCS open source finite element \sout{python} library is used to solve the variational forms \cite{AlnaesBlechta2015a}. The code is available at the following GitHub repository: \url{https://github.com/tpmatlab/THM-PFM}.

\subsection{Summary of Considered Cases}

Three sets of simulations were considered in the current work: i) first, a qualitative benchmark was proposed following the concrete column case proposed by Cheng et al. \cite{cheng2022coupled}; ii) second, the current methodology was used to investigate whether the thermomechanical stresses, which were shown to be critical for the case of fire accidents on Portland cement structures \cite{cheng2022coupled}, were also significant on the case of refractory castables drying which have much smaller heating rates; and, iii) the model was used to compare the effect of different heating methods. A list of the cases considered is presented in Table \ref{tab:THM-PF_experiments}.

\renewcommand{\arraystretch}{1.5}
\begin{table}[H]
\caption{Summary of the cases studied in the current work.}
\label{tab:THM-PF_experiments}
\centering
\resizebox{0.95\textwidth}{!}{%
\begin{tabularx}{\textwidth}{c c X X}
\toprule
\textbf{Case:} & \textbf{Model:} & \textbf{Description:} & \textbf{Objective:} \\
\toprule
\multirow{1}{*}[-2ex]{1} & \multirow{1}{*}[-2ex]{THM-PF} & Portland cement concrete column under fire. & Benchmark the current work implementation with Cheng et al. \cite{cheng2022coupled} \\
\midrule
2 & TH-PF & \multirow[c]{3}{=}[1.5ex]{Homogeneous heating of a refractory castable plane wall.} & \multirow[c]{3}{=}{Investigate the role of thermomechanical stresses in the explosive spalling during drying.} \\
3 & TM-PF & & \\
4 & THM-PF & & \\
\midrule
5 & THM-PF & {Single natural gas flame heating of a refractory castable plane wall.} & \multirow[c]{2}{=}[-4.5ex]{Compare different heating methods.} \\
6 & THM-PF & {Microwave heating of a refractory castable plane wall.} & \\
\bottomrule
\end{tabularx}
}
\end{table}

\begin{figure}[H]
\centering\includegraphics[width=0.85
\linewidth]{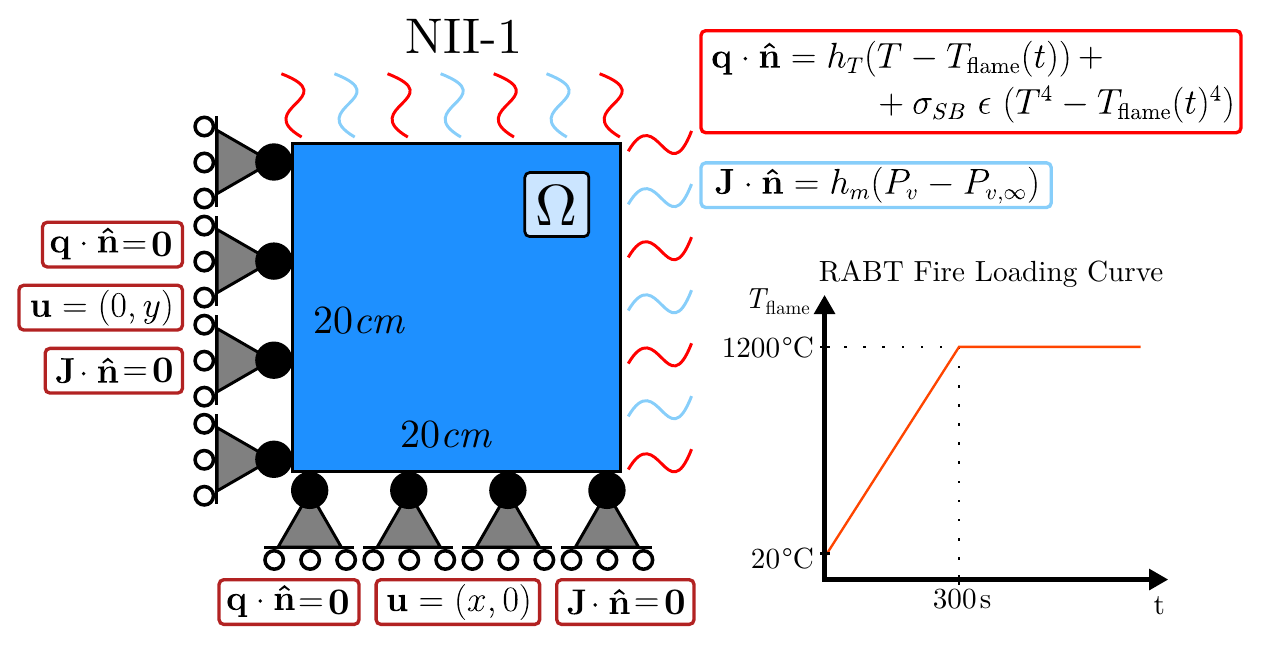}
\caption{Geometry and boundary conditions for the numerical benchmark test of the THM-PF model.} \label{fig:cheng_val_geo_BC}
\end{figure}

The configuration for the benchmark case is shown in Figure \ref{fig:cheng_val_geo_BC}, which depicts the behavior of a concrete column cross-section under the standard RABT heat-up curve, which assumes a linear temperature increase from 20°C to 1200°C over a 300 s period, after which the temperature is held at 1200°C \cite{cheng2022coupled}. 

The left and bottom faces are the symmetry planes of the column cross-section, whereas the right and top sides are permeable and are exposed to the RABT temperature increase through a Robin boundary condition. These two sides are also free to move (i.e., these surfaces have no displacement constraints). In order to simplify the formulation, a monophasic model is adopted for the thermohygro behavior, which is capable of predicting similar maximum pressures than the multiphasic models in previous studies by the authors \cite{moreira2021direct}.

\begin{figure}[H]
\centering
\includegraphics[width=1\linewidth]{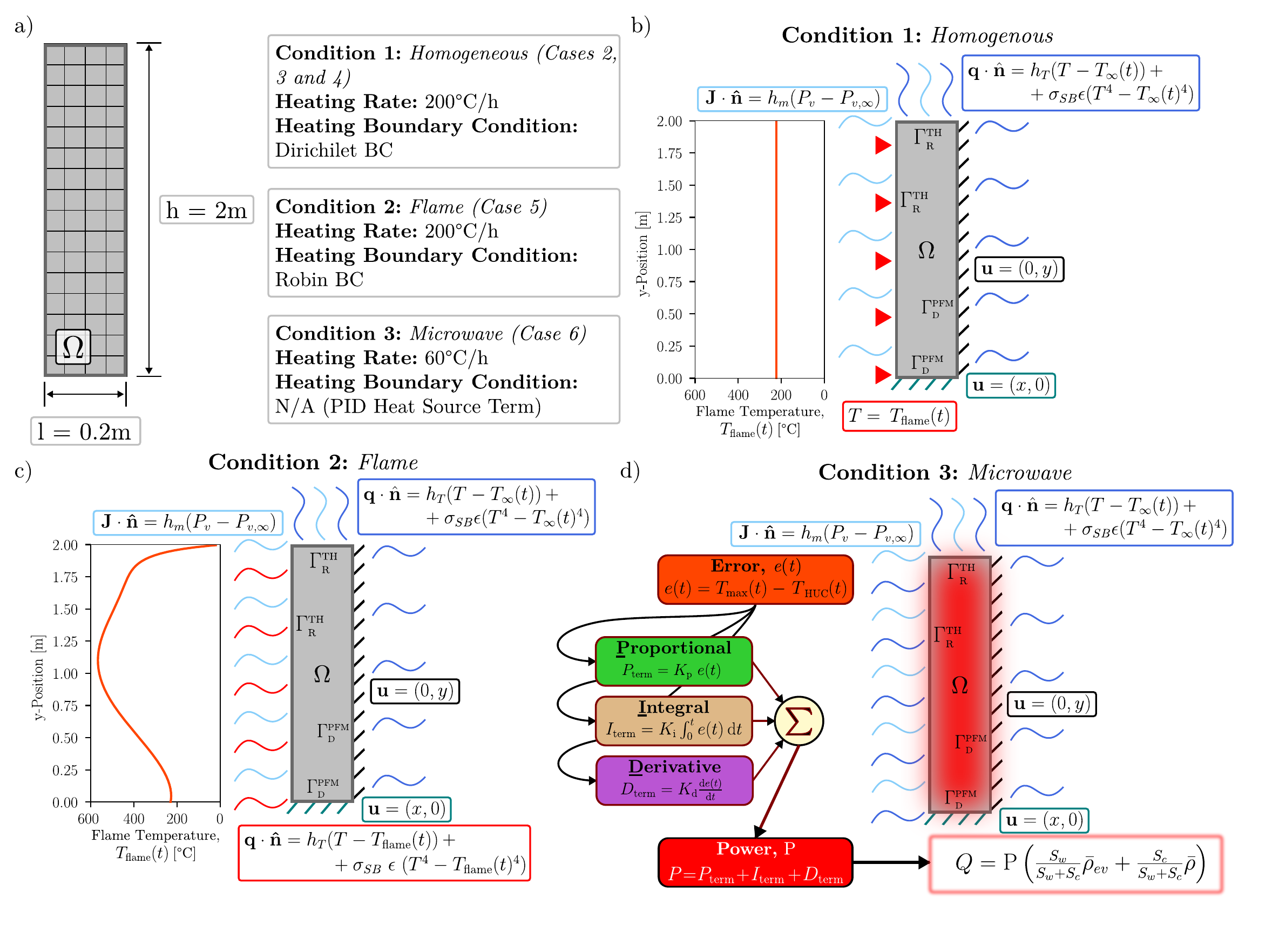}
\caption{a) Geometry and boundary conditions considered for the b) flame, c) homogeneous and d) microwave heating conditions. \label{fig:THMPF_cases}}
\end{figure}

The geometry and boundary conditions for tests 2 through 6 are summarized in Figure \ref{fig:THMPF_cases}. The coupling between the thermo and hygro parts of the model and the mechanical behavior were investigated in tests 2 and 3. For test 2, the thermomechanical stresses were disregarded, i.e. $\alpha = 0 \degree \textrm{C}^{-1}$, whereas for test 3, the vapor pressure effect was neglected by assuming a zero Biot coefficient ($b=0$).

The scenario of homogeneous heating is depicted in Figure \ref{fig:THMPF_cases} b), where a Dirichlet boundary condition is applied on the left face of a 2D section of the castable to represent the heating with a heating rate of 200°C/h. This face is also considered being permeable, as well as the boundary at the top of the lining, being free to transport mass to the surroundings. This top edge also loses thermal energy to the environment.

The left and right surfaces are also free to move. The boundary at the right represents the contact of the castable with a stiff steel shell (not represented in the simulations), which prevents horizontal displacement. This edge is also cooled by radiative heat transfer and natural convection Finally, the bottom of the lining is assumed to be adiabatic, impermeable and with the movement in the $y$-direction restricted.

The most common methodology adopted in the industry for drying refractory castables is based on gas burners \cite{Palmer2014Accelerated1}. This scenario is represented in Figure \ref{fig:THMPF_cases} c), where the temperature distribution applied follows the profiles reported by Liu et al. \cite{liu2021radiative} for a single oxy-fuel flame.

A Robin boundary condition represents the heating by radiative and natural convection. The decision to use a Robin boundary condition instead of a Dirichlet one (as applied on the homogeneous conditions) was made because the difference in temperature on the vertical direction would lead to regions with an effective heating rate higher than 200°C/h (the rate adopted in the homogeneous case). Therefore, by adopting the Robin boundary condition, the average heating rate for both conditions became comparable.

The last condition considered is the one depicted in Figure  \ref{fig:THMPF_cases} d), where the heating is provided by a microwave source. This is implemented following Gong et al. \cite{gong1991numerical}, by considering a volumetric heating source, with a variable input power $P$. The amount of energy absorbed is weighted by two factors,  one related to the microwave absorption by the evaporable water, $S_w$, and a second one, which  represents the refractory castable absorption. The proportionality constants are based on the values proposed by Gong et al. with $S_w = 1040$ [J/(Kg s)] and $S_c = 13$ [J/(Kg s)].

The input power is adjusted by a PID controller in order to match the maximum temperature within the castable, $T_{\textrm{max}}(t)$, with the corresponding temperature of the applied heat-up curve $T_{HUC}(t)$. To compensate the volumetric heating nature of the microwave approach, the heating rate applied  was decreased from 200°C/h (considered in cases 4 and 5) to 60°C/h. Finally, because the microwaves do not heat the surrounding air, in this case, the top, left and right faces are free to exchange heat with the surroundings at constant room temperature through natural convection and thermal radiation.

\section{Results and Discussions}
\label{sec:results}
\subsection{Benchmark of the THM-PF model}
\begin{figure}[H]
\centering\includegraphics[width=0.6\linewidth]{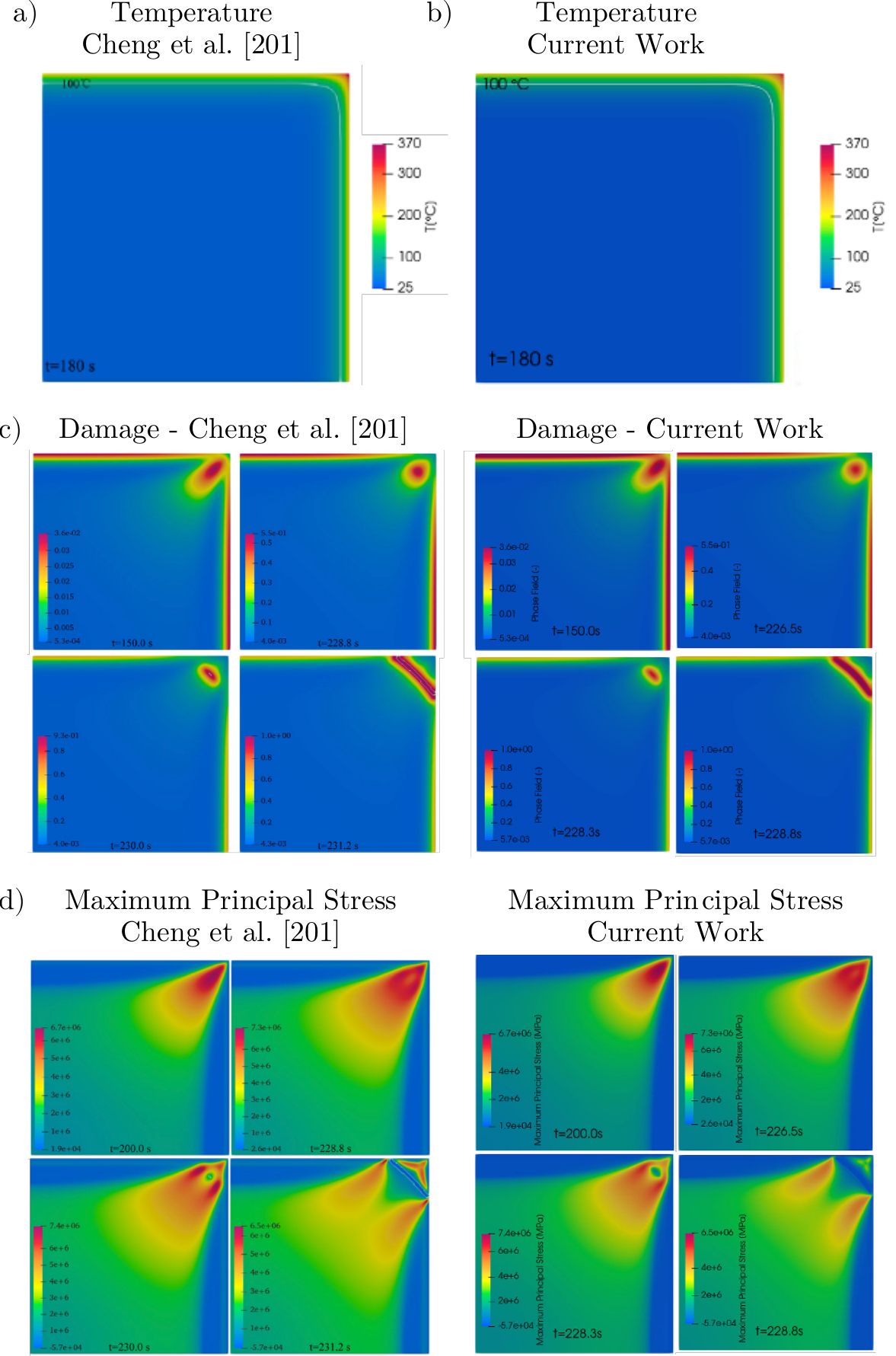}
\caption{Comparison of the results obtained by Cheng et al. \cite{cheng2022coupled} and those obtained by the anisotropic thermohygromechanical model implemented in the current work.} \label{fig:cheng_val}
\end{figure}

Considering the qualitative comparison with Cheng et al. \cite{cheng2022coupled}, the temperature distribution predicted by both models at the 180 seconds mark is shown in Figure \ref{fig:cheng_val} a) and b). Such results demonstrate that the thermal behavior calculated was the same in both models.

Figure \ref{fig:cheng_val} c) shows the damage evolution for both models, highlighting that in both cases the cracking took place in less than three seconds - between 228.8 and 231.2 seconds in the results by Cheng et al. and between 226.5 and 228.8 seconds in the current work. This is certainly a consequence of the assumption of quasi-static brittle behavior, as in the real cases, spalling is a violent process that occurs in fractions of a second \cite{Palmer2014Accelerated1, felicetti2024fast}.

\begin{figure}[H]
\centering\includegraphics[width=0.9\linewidth]{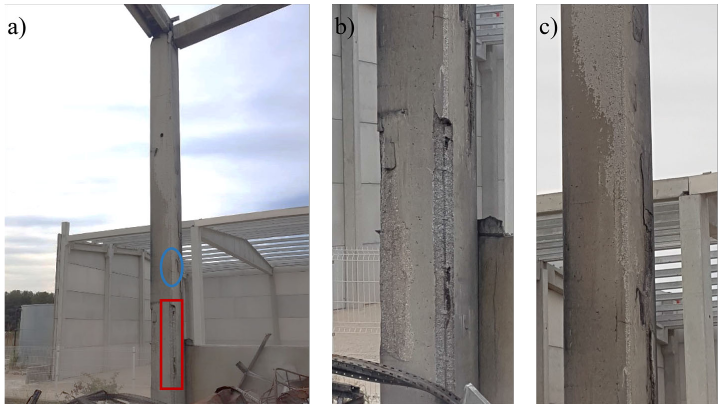}
\caption{a) Reinforced concrete column within an industrial structure post-fire; b) and c) show detailed views of corner spalling and surface deterioration on the RC column. From \cite{pena2023additional}.} \label{fig:col_exp}
\end{figure}

Nonetheless, both models predict similar cracking of the column corner thanks to the combination of thermal stresses and vapor pressurization, resulting in a similar failure mode as the one experimentally observed by et Peña et al. \cite{pena2023additional} and shown in Figure \ref{fig:col_exp}.

The maximum principal stresses shown in Figure \ref{fig:cheng_val} d) confirms the similarities between the results of both models. However, as with the damage case, there is a delay of 2.4 s, between the predictions, which can be explained by the different thermohygro model assumption of a multiphasic approach (as assumed by Cheng) and a single phase (adopted in this work).

Consequently, the present thermohygromechanical phase-field formulation is considered correctly implemented, as it reproduces the behavior reported in experimental tests and in previous models.

The main difference that highlights the advantages of the current implementation is on the total simulation time, described in Table \ref{tab:cpu_time}. The comparison reveals a substantial reduction in computational cost when adopting the adaptive time-stepping strategy proposed in this work.

While the fixed time-step scheme required 4 h 42 min to complete the simulation, the adaptive algorithm reduced this time to 50 min (an improvement in computational cost by a factor of $\sim 5.7$). 

\begin{table}[H]
\centering
\caption{Comparison of total simulation times obtained using fixed and adaptive time step algorithms.}
\label{tab:cpu_time}
\begin{tabular}{lc}
\toprule
                   & \multicolumn{1}{l}{Simulation Time} \\
\toprule
Fixed Time Step    & 4 h 42 min                          \\
Adaptive Time Step & 50 min     \\
\bottomrule
\end{tabular}
\end{table}

This gain reflects the algorithm’s ability to refine the time step only when required by the evolving thermohygromechanical fields, thereby avoiding unnecessary evaluations during stages where minor changes occur, which is critical for simulating the long drying processes of refractory castables. 

\subsection{The Role of Thermomechanical Stresses on the Drying of Castables}

\begin{figure}[H]
\centering\includegraphics[width=0.85\linewidth]{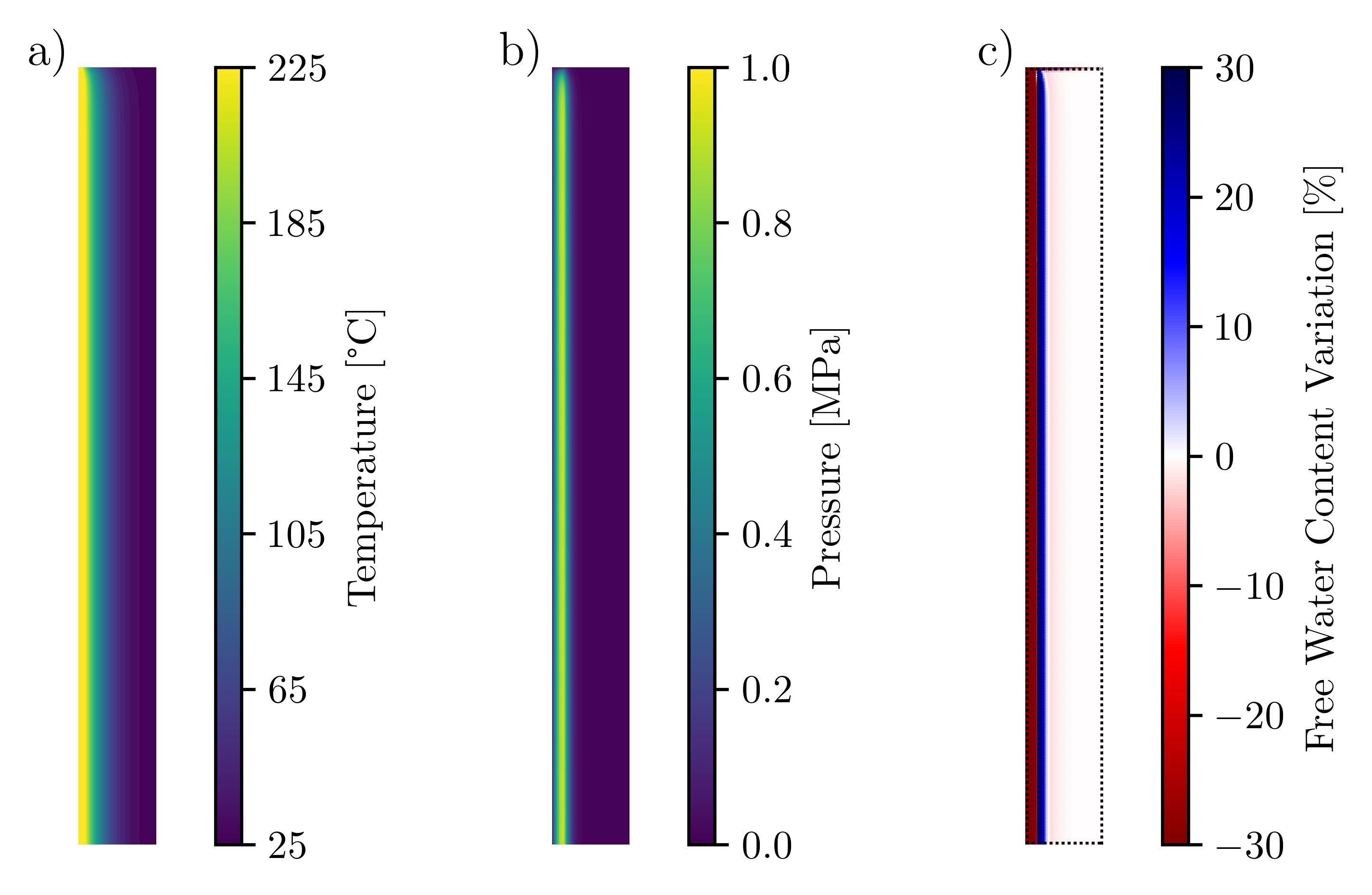}
\caption{a) Temperature, b) moisture pressure and c) free water content distributions at t = 85.8 min, the moment when there was crack initiation in Case 4.} \label{fig:TM_effect_PT}
\end{figure}

With respect to the role of thermomechanical stresses on the drying of castables, Figure \ref{fig:TM_effect_PT} shows the distribution of temperature, pressure and free water at the moment of crack initiation as predicted in Case 4 (the full THM-PF model). As the current implementation considers a weak coupling between the mechanical and thermohygral behavior (i.e. the thermohygral problem is not impacted by the mechanical behavior), the case where the thermomechanical stresses are disregarded (Case 2) and the case where the Biot's coefficient is assumed to be zero (Case 3) also share the same distribution of temperature, pressure and free water content at the moment of crack initiation. 

The hot face reached a temperature of 311.4 °C, at the 85.8 min mark, as shown in Figure \ref{fig:TM_effect_PT} a). At this moment, the pressure front had penetrated roughly 3 cm inside the castable wall, with a maximum value of 1 MPa, as shown on Figure \ref{fig:TM_effect_PT} b). The moisture exchange with the environment taking place at the top surface of the wall made the pressure front to be slightly shifted towards the bottom of the domain.

This drying front had concentrated water on the regions ahead of it, following the pattern expected by the moisture clogging phenomenon \cite{moreira2022experimental}, as depicted on Figure \ref{fig:TM_effect_PT} c).

\subsection{Influence of Thermohygral Strains}

\begin{figure}[H]
\centering\includegraphics[width=0.9\linewidth]{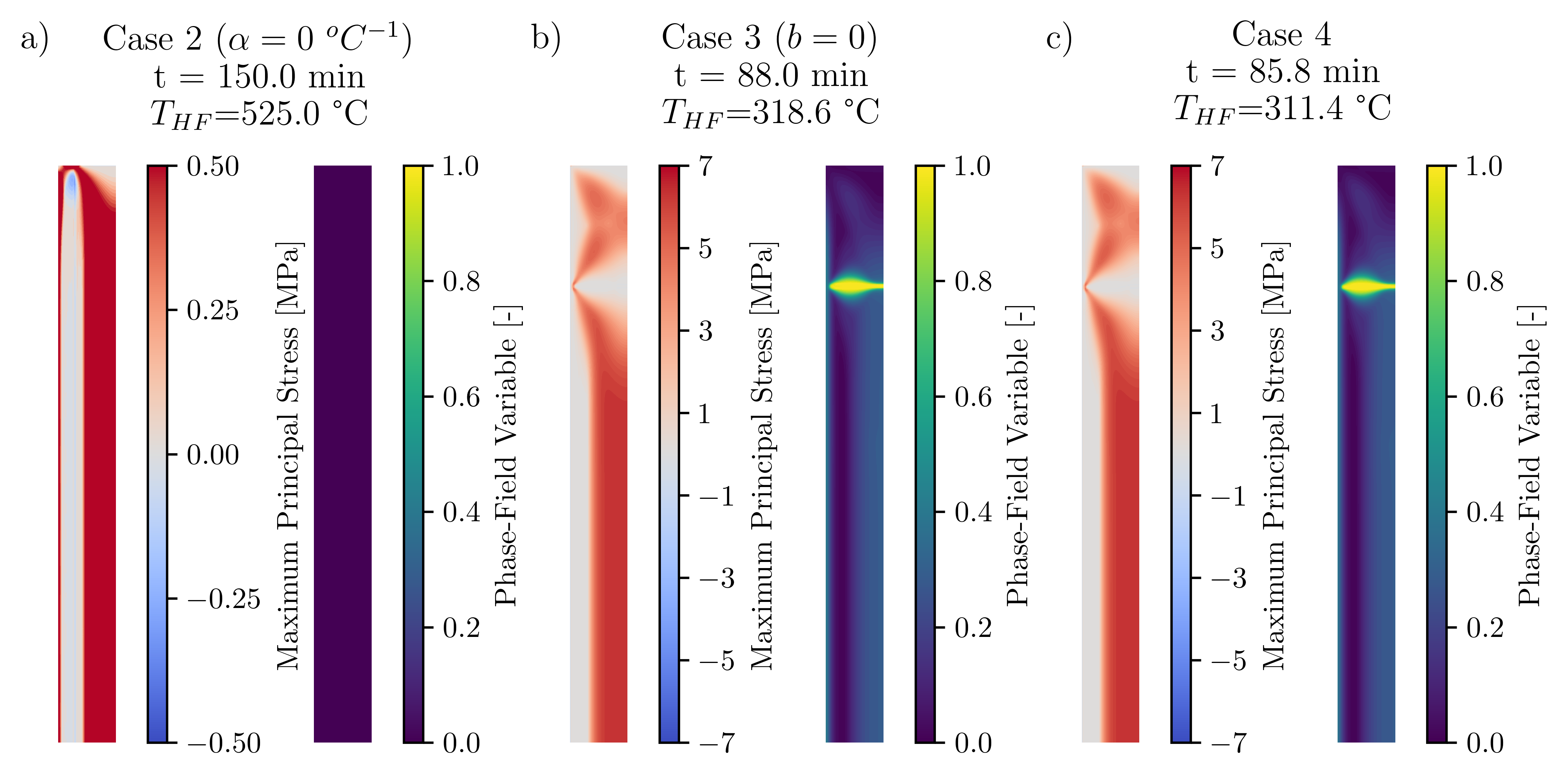}
\caption{Distribution of the principal maximum stresses and the damage variable, $\phi$, at the moment the material failed when considering a) only thermohygral strains, b) thermomechanical ones, and c) both.} \label{fig:TM_effect_cracks}
\end{figure}

Figure \ref{fig:TM_effect_cracks} highlights the maximum principal stress and the phase field variable for all the cases. On case 2, the mechanical stress is driven exclusively by the internal vapor pressure gradients generated during drying. The maximum principal stress reveals a highly localized tensile zone ($\approx 0.5\text{ MPa}$) along the vertical sides. This is a consequence of the localized internal pressure which tries to expand the castable wall, inducing tensile stresses in the solid skeleton. There is also a small compressive region at the center of the top section ($\approx -0.1\text{ MPa}$), resulting from the expansion of the material surrounding this position.

Due to the relatively small values of maximum principal strain shown in Figure \ref{fig:TM_effect_cracks} a), no cracks appeared during the simulation, in contrast with the results of the thermomechanical and thermohygromechanical analysis depicted in Figures \ref{fig:TM_effect_cracks} b) and \ref{fig:TM_effect_cracks} c). This indicates that while thermohygral strains contribute to the overall stress state, the pore pressure alone, considering the assumed boundary conditions and geometry, is insufficient to trigger cracking. 

The contribution of thermomechanical strains increased the resulting maximum principal stress by an order of magnitude, reaching values of approximately $7.0\text{ MPa}$ in both cases 3 and 4. This significant increase in stress levels is attributed to the thermal expansion of the material, which, when combined with the mechanical constraints imposed by the boundary conditions, leads to a pronounced internal bending moment and subsequent crack initiation.

This bending moment is a direct consequence of the displacement restrictions on the bottom and right-hand side of the domain, which prevent free thermal expansion in the lower part of the wall. As the left upper portion remains unconstrained it can expand freely. This kinematic asymmetry results in a severe internal bending moment as the exposed left face undergoes thermal expansion. Consequently, a massive concentration of tensile stress develops precisely at the transition zone ($y \approx 1.3\text{ m}$), resulting in a localized horizontal macro-crack.

It is worth noting that in industrial practice the explosive spalling typically takes place as cracks running parallel to the heated surface rather than the horizontal macro-cracking predicted here. This difference exists because actual industrial linings are subjected to severe lateral and vertical confinement from neighboring refractories. This induces massive compressive stresses parallel to the hot face, effectively preventing the initiation of horizontal cracks and forcing the material to fail by buckling or delamination. 

Nevertheless, this does not invalidate the comparisons between the cases. By adopting an unconfined configuration, the model can better isolate the effects of the different mechanisms, which could be masked by the hydrostatic confinement on real refractory linings. Thus, the current results highlight the relative importance of the thermomechanical strains in the drying of refractory castables both in the simplified analysis and in more realistic scenarios.

These results also agree with the current understanding that explosive spalling arises from the combined action of thermal stresses, vapor pressurization, and potentially additional mechanisms such as flash vaporization \cite{Ozawa2012Study, shen2020multiphysics, cheng2022coupled, felicetti2024fast}. These mechanisms were originally proposed to explain spalling in Portland cement-based structures exposed to fire, and even though in such cases the heating rate is significantly higher, they also appear to govern the behavior of refractory castables during drying failure.

\subsection{The Comparison Between Different Heating Methods}
\begin{figure}[!htpb]
\centering\includegraphics[width=0.8\linewidth]{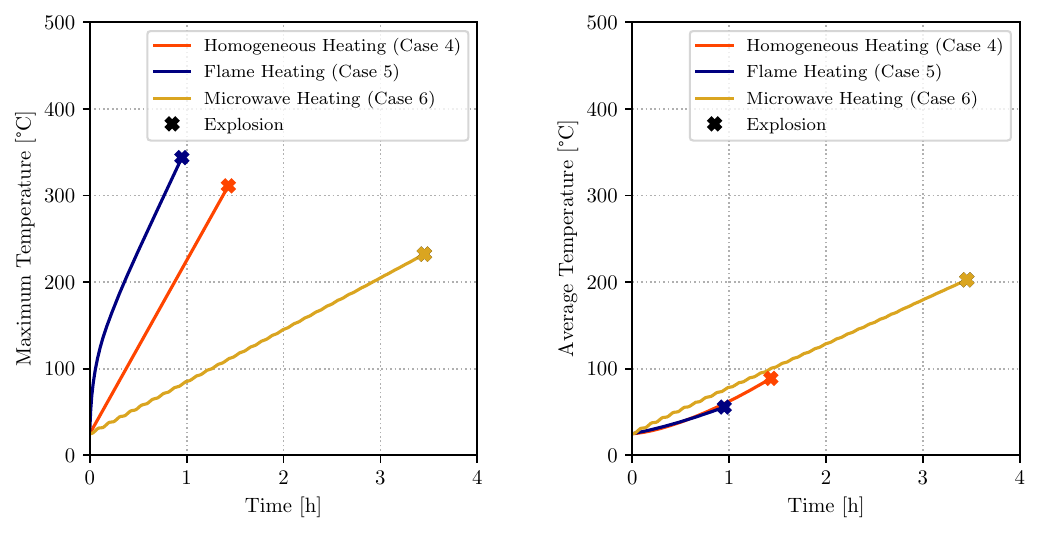}
\caption{Temporal evolution of a) maximum temperature in the domain and b) average temperature of the domain.} \label{fig:PF-temps}
\end{figure}

The choice of the right metric to compare different heating methods can be a complex task, as  illustrated in Figure \ref{fig:PF-temps}. If the evolution of the maximum temperature on the wall is considered, the scenario of flame heating could be viewed as the best alternative, as it could reach the highest temperature before the moment of failure, as shown in Figure \ref{fig:PF-temps} a). This occurs because the non-uniform heating has a maximum temeperature that is larger than the values observed in the uniform heating scenario. Additionally, the microwave heating case attained the lowest maximum temperature, thanks to its lower heating rate (60°C/h vs 200°C/h of the other methods).

Nevertheless, due to the fact that the microwave heating is volumetric, the measure that is more representative is the evolution of the average temperature of the domain, as shown in Figure \ref{fig:PF-temps} b). In this case, it is clear that the adoption of the method with the smallest induced thermal stresses (i.e. the microwave heating) allowed for a higher average temperature to be reached before failure, even though it was the method with the slowest heating rate.

\begin{figure}[H]
\centering\includegraphics[width=0.8\linewidth]{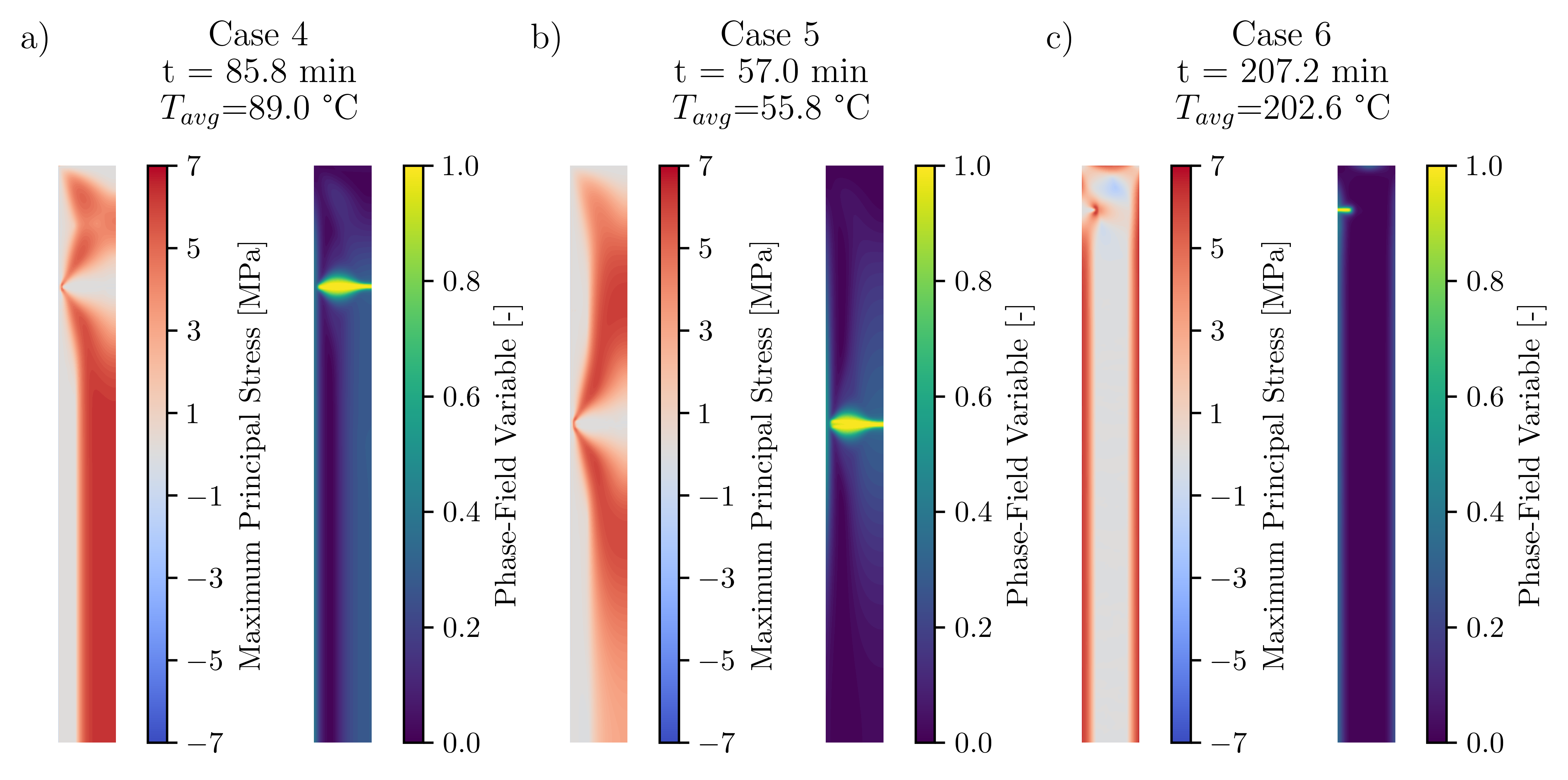}
\caption{Distribution of the principal maximum stresses and the damage variable, $\phi$, at the moment the material failed under heating conditions of a) homogeneous, b) flame, and c) microwave.} \label{fig:PF-Cracks}
\end{figure}

Next, the distribution of the principal maximum stresses along the wall, as well as the phase-field variable at the moment of failure is shown for the different heating scenarios in Figure \ref{fig:PF-Cracks}. The overall distribution of principal stress for the flame heating scenario is similar to the one observed for the homogeneous case, with the difference being that the material fails at the region of maximum temperature (i.e. the middle of the wall) and at an earlier moment.

On the other hand, the volumetric heating of the microwave case combined with the convective cooling of the lining surfaces resulted in opposite thermal gradients on the wall, which resulted in tensile stresses on its surface, leading to the formation of cracks. In an industrial setting, such challange can be overcome by providing additional heating through hot dry air, which also enhances the water withdrawal, as described by Goto et al. \cite{goto2013development}.

\begin{figure}[!htpb]
\centering\includegraphics[width=0.6 \linewidth]{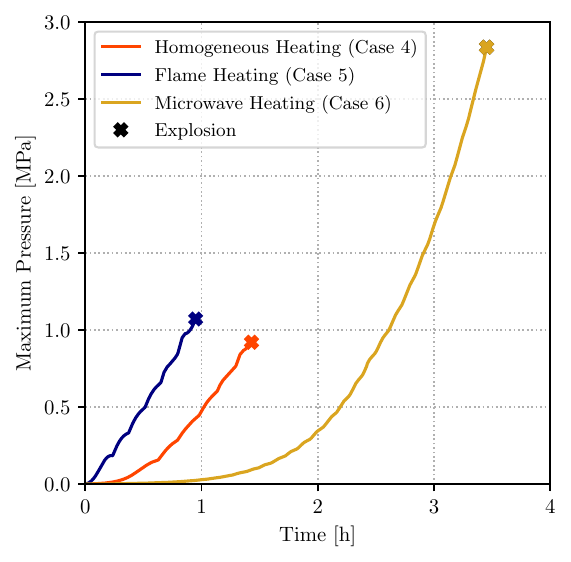}
\caption{Temporal evolution of the maximum pressure calculated in the different heating scenarios.} \label{fig:PF-max_pressure}
\end{figure}

The improved resistance against drying failure of the microwave heating scenario is also highlighted by the evolution of the maximum pressure within the castable wall for all the heating cases, as shown in Figure \ref{fig:PF-max_pressure}. The homogeneous heating scenario resulted in a maximum pressure of 0.9 MPa at the moment of failure, whereas the flame heating resulted in 1.1 MPa. The microwave heating scenario, on the other hand, could withstand pressure levels as high as 2.8 MPa before failure.

Thus, through the analysis of the results in the current configuration, it is possible to conclude that the usually adopted procedure of heating refractory castables with a single flame directly reduces the effective refractory castable resistance to pressurization, which can lead to explosive spalling.

\clearpage

\section{Conclusions}
\label{sec:conclusion}

The current work presented an open-source implementation of a thermohygromechanical phase-field model for studying the high temperature behavior of refractory castables and Portland cement concrete structures. The model is readily available for download and use, and it is based on the FEniCS open-source finite element library. 

The numerical model was based on the assumption of quasi-static brittle fracture, with a hybrid monolithic-staggered solution strategy, which enabled the use of an adaptive pseudo-time step strategy capable of reducing the total simulation time by a factor of five for the studied cases. The results of the model were also qualitatively validated against both experimental observations and numerical results from the literature.

Additionally, the analysis of refractory castable drying using this tool revealed that even in this case, where the heating rates are significantly lower than those observed in fire accidents, the thermomechanical stresses play a non-negligible role in the explosive spalling behavior. 

Finally, the comparison between different heating methods available in the industry, namely, homogeneous, flame, and microwave heating, have demonstrated that alternatives that minimize thermal gradients inside the material are promising to mitigate the risk of failure during drying, improving the overall safety and efficiency of refractory castables installation without the need to change the material composition.

\section{Acknowledgments}\label{sec:acknowledgments}
This study was financed in part by the Coordenação de Aperfeiçoamento de Pessoal de Nível Superior - Brasil (CAPES) - Finance Code 001 and the ANR (Agence Nationale pour la Recherche) MULTI-FIRE n.ANR-23-CE51-0001. This study was financed, in part, by the São Paulo Research Foundation (FAPESP), Brasil (grants 2021/00251-0 and 2025/00975-9). R.F.A. is partially funded by National Council for Scientific and Technological Development (CNPq, Brazil), grant 308704/2022-3.

\section{Data Availability}\label{sec:data_availability}
The data that support the findings of this study are available at the GitHub repository: \url{https://github.com/tpmatlab/THM-PFM}.

\section{CRediT Authorship Contribution Statement}\label{sec:credit}
\noindent\textbf{Murilo Henrique Moreira:} Conceptualization, Methodology, Formal Analysis, Software, Data curation, Writing – original draft, review \& editing. \\
\textbf{Roberto Federico Ausas:} Methodology, Software, Investigation, Data curation, Visualization, Writing – review \& editing.\\
\textbf{Stefano Dal Pont:} Funding acquisition, Supervision, Conceptualization, Methodology, Supervision, Writing – review \& editing. \\






\clearpage
\bibliographystyle{elsarticle-num-names}
\bibliography{references.bib}






\doublespacing

\clearpage
\appendix

\section{Properties and Constitutive Equations of the 5CAC} \label{sec:app_const_eq_5CAC}

\subsection*{Elastic Properties}

The elastic modulus $E$ and critical energy release rate $G_c$ are assumed constant:
\begin{equation}\label{eq:5CAC_E}
    E = E_0 = 14 \times 10^{9} \text{ Pa}
\end{equation}

\begin{equation}\label{eq:5CAC_Gc}
    G_c = G_{c,0} = 95 \text{ N/m}
\end{equation}

The Poisson's ratio is $\nu = 0.25$. The Lamé constants for the plane strain assumption are:
\begin{equation}\label{eq:5CAC_mu}
    \mu = \frac{E}{2(1+\nu)}
\end{equation}
\begin{equation}\label{eq:5CAC_lambda}
    \lambda = \frac{E \nu}{(1+\nu)(1-2\nu)}
\end{equation}

The bulk modulus is:
\begin{equation}\label{eq:5CAC_K}
    K = \frac{3\lambda + 2\mu}{3}
\end{equation}

The thermal expansion coefficient is $\alpha_T = 1 \times 10^{-5}$ K$^{-1}$ and the characteristic length scale is $\ell_0 = 4.1 \times 10^{-3}$ m.

\subsection*{Thermohygral Properties}

The saturation vapor pressure $P_s$ (Pa) is given by the empirical formula \cite{bazant1978}:
\begin{equation}\label{eq:5CAC_Ps}
    P_s(T) = \exp\!\left(\frac{C_1}{T} + C_2 + C_3 T + C_4 T^2 + C_5 T^3 + C_6 \ln T\right)
\end{equation}
where $C_1 = -5800.2206$, $C_2 = 1.3914993$, $C_3 = -4.8640239 \times 10^{-2}$, $C_4 = 4.1764768 \times 10^{-5}$, $C_5 = -1.4452093 \times 10^{-8}$, and $C_6 = 6.5459673$, and $T$ is the absolute temperature in Kelvin.

The relative humidity is accounted for as:
\begin{equation}\label{eq:5CAC_h}
    h(P_v, T) = \frac{P_v}{P_s(T)}
\end{equation}

The dehydration water content $w_d$ (kg/m$^3$) as a function of temperature (in °C, $T_{C} = T - 273.15$) is represented by a sigmoid function:
\begin{equation}\label{eq:5CAC_wd}
    w_d(T) = A_1 + \frac{A_2 - A_1}{1 + \exp\!\left(\dfrac{T_C - x_0}{d_x}\right)}
\end{equation}
where $A_1 = 22.839$ kg/m$^3$, $A_2 = 0.139$ kg/m$^3$, $x_0 = 271.39$°C, and $d_x = 23.06$°C.

The sorption isotherm describes the evaporable water content $w$ (kg/m$^3$). The parameter $m(T)$ governing the sorption curve shape is:
\begin{equation}\label{eq:5CAC_m}
    m(T) = 1.04 - \frac{(T - 263.15)^2}{(T - 263.15)^2 + 22.34\,(298.15 - 263.15)^2}
\end{equation}

The sorption isotherm is defined by two piecewise parts. For $h \leq 0.96$:
\begin{equation}\label{eq:5CAC_w1}
    w_1(P_v, T) = w_c \left(\frac{w_0}{w_c}\, h(P_v, T)\right)^{1/m(T)}
\end{equation}

For $h > 1.04$:
\begin{equation}\label{eq:5CAC_w2}
    w_2(P_v, T) = w_c \left[0.037\,(h - 1.04) + 0.3335\left(1 - \frac{T_C^2}{3.6 \times 10^5}\right)\right]
\end{equation}

In the transition region $0.96 \leq h \leq 1.04$, a linear interpolation between $w_1(0.96\,P_s, T)$ and $w_2(1.04\,P_s, T)$ is used. The full sorption isotherm is then:
\begin{equation}\label{eq:5CAC_w}
    w(P_v, T) =
    \begin{dcases}
        w_1(P_v, T), & h \leq 0.96 \\
        w_1(0.96\,P_s, T) + \dfrac{\bigl[w_2(1.04\,P_s, T) - w_1(0.96\,P_s, T)\bigr](h - 0.96)}{0.08}, & 0.96 < h \leq 1.04 \\
        w_2(P_v, T), & h > 1.04 \\
        0, & T \geq 647.15 \text{ K}
    \end{dcases}
\end{equation}
where $w_c = 300$ kg/m$^3$ is the cement content, $w_0 = 100$ kg/m$^3$ is the initial water content, and $T_C = T - 273.15$.

The solid density is $\rho_s = 2200$ kg/m$^3$, the specific heat of the solid is $C_p = 1100$ J/(kg K), and the specific heat of liquid water is $C_{p,l} = 4100$ J/(kg K). The effective thermal conductivity is $\lambda_c = 1.67$ W/(m K).

The latent heat of evaporation is given by:
\begin{equation}\label{eq:5CAC_DeltaHev}
    \Delta H_{ev}(T) =
    \begin{dcases}
        3.5 \times 10^{5}\,(647.3 - T)^{1/3}, & T \leq 647.3 \text{ K} \\
        0, & T > 647.3 \text{ K}
    \end{dcases}
\end{equation}

The hydraulic conductivity $a(P_v, T)$ (m/s) follows the formulation of Bažant \cite{bazant1978}. For $T \leq 368.15$ K:
\begin{equation}\label{eq:5CAC_a_low}
    a(P_v, T) = a_0\, f_1(T, h)\, f_2(T)
\end{equation}
and for $T > 368.15$ K:
\begin{equation}\label{eq:5CAC_a_high}
    a(P_v, T) = 5.6\, a_0\, f_3(T)
\end{equation}
where $a_0 = 10^{-14}$ m/s is the initial hydraulic conductivity. The humidity-dependent part is:
\begin{equation}\label{eq:5CAC_f1}
    f_1(T, h) =
    \begin{dcases}
        a_t(T) + \dfrac{1 - a_t(T)}{1 + \bigl[4(1-h)\bigr]^4}, & h < 1 \\
        1, & h \geq 1
    \end{dcases}
\end{equation}
with:
\begin{equation}\label{eq:5CAC_at}
    a_t(T) = \frac{0.95\,T_C}{70} - 0.28929
\end{equation}

The temperature-dependent coefficients are:
\begin{equation}\label{eq:5CAC_f2}
    f_2(T) = \exp\!\left[2700\left(\frac{1}{298.15} - \frac{1}{T}\right)\right]
\end{equation}
\begin{equation}\label{eq:5CAC_f3}
    f_3(T) = \exp\!\left(\frac{T_C - 95}{0.881 + 0.214\,(T_C - 95)}\right)
\end{equation}

\end{document}